\documentclass[12pt]{amsart}
\usepackage{amsmath,amssymb,amsfonts,amscd,graphicx,xypic}

\newtheorem{lemma}[subsection]{Lemma}
\newtheorem{theorem}[subsection]{Theorem}

\theoremstyle{definition}

\theoremstyle{remark}
\numberwithin{equation}{section}
\theoremstyle{plain}
\newtheorem{thm}{Theorem}

\theoremstyle{definition}

\begin{document}

\title{Sticky Cantor Sets in ${\mathbf{\mathbb R}^d}$}

\author{Vyacheslav Krushkal}

\address{Department of Mathematics, University of Virginia,
Charlottesville, VA 22904-4137 USA}
\email{krushkal\char 64 virginia.edu}

\begin{abstract} 
A subset of ${\mathbb R}^d$ is called ``sticky'' if it cannot be isotoped off of itself by a small ambient isotopy. Sticky wild Cantor sets are constructed in ${\mathbb R}^d$ for each $d\geq 4$.
\end{abstract}

\maketitle

\section{Introduction}

Wild  Cantor sets in Euclidean spaces have been extensively studied, following the definition of Antoine's necklace in \cite{A}. Embeddings exhibiting various interesting phenomena have been constructed using decomposition theory, starting with the Bing decomposition of the $3$-sphere \cite{B}. The reader may find a general discussion and examples in \cite{D, DV}.

A systematic construction of Cantor set embeddings which suitably approximate codimension $2$ submanifolds was given in \cite{DE}. Accordingly, a notion of ``general position'' for wild Cantor sets in Euclidean spaces may seem plausible (and has been mentioned in the literature). Indeed, the result of this paper may be thought of as a version of general position for certain wild Cantor sets approximating codimension $2$ submanifolds. (It may be interesting to compare this with \cite{CW} where wild Cantor sets are constructed in ${\mathbb R}^d$, $d\geq 3$, which can be ``slipped off'' every Cantor set.)

This paper concerns the question of  whether there exist  sticky Cantor sets in ${\mathbb R}^d$  (see \cite[Conjecture 1]{W}, and also \cite[Problem E8]{D1}). Here a (wild) Cantor set embedded in ${\mathbb R}^d$ is called {\em sticky}  if it cannot be isotoped off of itself by any sufficiently small ambient isotopy. 
The question of whether any given Cantor set $X$ can be slipped off {\em every} Cantor set $Y$ in ${\mathbb R}^d$ has also been asked in \cite[Conjecture 1]{D0}, \cite[Conjecture 1.2]{CW}. The following theorem states the main result of the paper.

\begin{thm} \label{theorem} \sl 
For any $d\geq 4$ there exists a sticky Cantor set in ${\mathbb R}^d$.
\end{thm}

The main ingredients of the proof are the spun Bing decompositions, considered in \cite{Edwards, L}, and the Stallings theorem \cite{S} on the lower central series of groups. The spun Bing decompositions produce wild Cantor sets which approximate codimension $2$ submanifolds (more specifically, $S^{d-2}\subset {\mathbb R}^d$), in the sense of \cite{DE}. It is an interesting question whether a similar ``general position'' result holds for the construction in \cite{DE}.

\section{Proof of theorem \ref{theorem}}

\begin{figure}[ht]
\includegraphics[height=3.5cm]{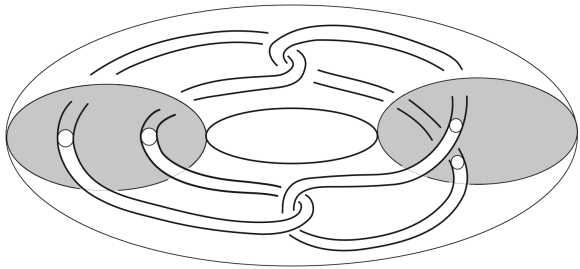}
 \caption{Bing double of the core of a solid torus. The spun version is obtained by spinning this figure in ${\mathbb R}^4$ about the $2$-plane which intersects the solid torus in the two indicated shaded disks.}
\label{Bing figure}
\end{figure}

{\em Proof.} First consider the case $d=4$.
The starting point of the construction is the nested sequence of Bing doubles \cite{B} in the solid torus, figure \ref{Bing figure}.
Spin it to obtain a nested sequence of ``Bing doubles of spheres'' in $S^2\times D^2$. (Each stage in this sequence is a collection of $S^2\times D^2$'s.) A well-known theorem \cite{Edwards}, \cite{L} (used in particular in the proof of the Double Suspension Theorem in \cite{Edwards}) states that  the spun Bing decomposition is shrinkable, so the nested sequence gives rise to a Cantor set in $S^2\times D^2$. In fact, it is shown in \cite{Edwards} that iterated spun Bing decompositions are shrinkable, yielding a Cantor set in $S^n\times D^2$, for each $n\geq 2$. 

\smallskip

\begin{figure}[ht]
\includegraphics[height=3.1cm]{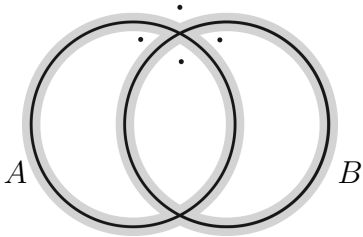}
\put(-127,20){$A$}
\put(-1,20){$B$}
 \caption{Two intersecting spheres in ${\mathbb R}^4$ and the Clifford torus.}
\label{spheres}
\end{figure}

Now consider two standard $2$-spheres in ${\mathbb R}^4$ intersecting in two points, figure \ref{spheres}, and let $A, B$ denote two copies of the Cantor set 
as above, one in each $S^2\times D^2$.  The proof of theorem \ref{theorem} relies on the following lemma.

\smallskip

\begin{lemma} \label{lemma} \sl
  $A$ cannot be isotoped off of $B$ by any sufficiently small ambient isotopy.
  \end{lemma}

 {\em Proof of lemma \ref{lemma}.} 
Consider a Clifford torus for one of the intersection points of the spheres. (The intersection of spheres is locally modeled on ${\mathbb R}^2\times \{ 0\}\, \cup\,  \{ 0\} \times {\mathbb R}^2 \, \subset \, {\mathbb R}^4$. The Clifford torus of radius $r$ is the product $S^1_r\times S^1_r$ of two circles of radius $r$  in the complement of these two planes in ${\mathbb R}^4$. The radius $r$ is chosen so that the Clifford torus is in the complement of the two $S^2 \times D^2$'s, so it is in the complement of $A\cup B$. The  torus is drawn as $4$ points ($S^0\times S^0$) in the $2$-dimensional illustration in figure~\ref{spheres}.)

Consider the two standard generators of ${\pi}_1$(Clifford torus): meridians $m_A$, $m_B$ to the two $2$-spheres, in other words the circles $S^1_r\times *$, $*\times S^1_r$ linking the $2$-spheres.
The $2$-cell of the torus gives the relation $[m_A,m_B]=1$. 

Suppose an isotopy in the statement of lemma 2 exists. Then some finite stages $A_n, B_n$ of the nested sequences defining $A, B$ are disjoint after the isotopy. The ambient isotopy is assumed to be small enough so that the nested sequences defining $A, B$ stay disjoint from the Clifford torus during the isotopy. To avoid confusion, we will keep the notation $A_n, B_n$ for the $n$-th stages of the nested sequences {\em before} the isotopy, and $\widetilde A_n$ will denote $A_n$ {\em after} the isotopy. The assumption is $\widetilde A_n\cap B_n=\emptyset$.

It follows from the previous two paragraphs that 
\begin{equation} \label{relation}
[m_A,m_B]=1 \; \, {\rm in} \; \, {\pi}_1({\mathbb R}^4\smallsetminus(\widetilde A_n\sqcup B_n)).
\end{equation}

We will next show that this leads to a contradiction with the Stallings theorem.  (The argument below is a generalization of the analogous proof that the two $2$-spheres shown in figure \ref{spheres} (as opposed to their ``aproximations'' $A_n, B_n$) cannot be made disjoint by a small isotopy. The reader may consider applying the argument in this case before carrying it our in full generality, discussed below.)

Since $\widetilde A_n\sqcup B_n$ is a collection of disjointly embedded (thickened) $2$-spheres, by Alexander duality $H_1({\mathbb R}^4\smallsetminus(\widetilde A_n\sqcup B_n) ; {\mathbb Z})$ is generated by meridians (small linking circles) to the $2$-spheres, and $H_2({\mathbb R}^4\smallsetminus(\widetilde A_n\sqcup B_n) ; {\mathbb Z})$ is trivial.
Given a group $G$, its lower central series is defined inductively by $G^1= G, \, G^2 = [G, G], \, . . . , \, G^k = [G, G^{k-1}]$.

\begin{theorem}[Stallings Theorem \cite{S}] \sl
 Suppose a map $f\! : X\longrightarrow Y$  induces
an isomorphism on $H_1(. ; {\mathbb Z})$ and a surjection on  $H_2(. ; {\mathbb Z})$. Then for each finite $k$,
$f$ induces an isomorphism ${\pi}_1(X)/({\pi}_1(X))^k \cong {\pi}_1(Y )/({\pi}_1(Y ))^k$.
\end{theorem}

For brevity of notation denote $Y:={\mathbb R}^4\smallsetminus(\widetilde A_n\sqcup B_n)$.
Consider the map $f\! \! :\vee^N S^1\longrightarrow Y$ from the wedge of circles (one circle for each $2$-sphere in the collection $\widetilde A_n\sqcup B_n$), mapping each circle to a meridian (small linking circle, connected to a basepoint) of the corresponding $2$-sphere. 
 This map induces an isomorphism on $H_1$ and a surjection on $H_2=0$. 
By the Stallings theorem, for any $k$,
$ {\pi}_1(Y )/({\pi}_1(Y ))^k$ is isomorphic to the corresponding quotient of the free group, $ F_N/(F_N)^k$ where $F_N\cong {\pi}_1(\vee^N S^1)$ is the free group on $N$ generators.

\begin{figure}[ht]
\includegraphics[width=4.5cm]{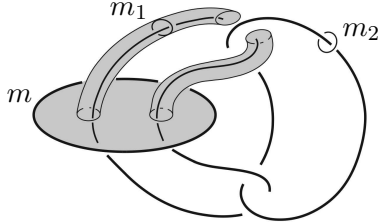}
{\small
\put(-11,71){$m_2$}
\put(-99,79){$m_1$}
\put(-138,46){$m$}
}
\caption{$m=[m_1,m_2]$, where $m$, $m_1$, $m_2$ are meridians suitably connected to a basepoint.}
\label{calc figure}
\end{figure}

The meridians $m_A, m_B$ may be explicitly written down both in ${\pi}_1({\mathbb R}^4\smallsetminus(A_n\sqcup B_n))$ and in ${\pi}_1(Y)$.
Indeed, each Bing doubling replaces a meridian $m$ to a sphere by a commutator of the meridians $m_1, m_2$ of the two smaller spheres. This may be read off in the $3$-space slice (shown in figure \ref{Bing figure}), as illustrated in figure \ref{calc figure}: $m$ bounds a punctured torus in the complement of the two Bing-doubled curves, with the two generators of ${\pi}_1$ of the torus corresponding to $m_1, m_2$. (This is an important basic example in Milnor's theory of link homotopy \cite{M}.)

The finite stages $A_n, B_n$ of the nested sequences defining $A, B$ correspond to a certain number of iterations of Bing doubling. This exhibits the meridians $m_A, m_B$ as commutators of fixed length $l$ in the meridians to the $2$-spheres forming $A_n, B_n$, where $l$ is the number of $2$-spheres in the each of the collections $A_n, B_n$.  (To relate this to the notation following the statement of the Stallings theorem, $N=2l$.)

\begin{figure}[ht]
\includegraphics[width=5.5cm]{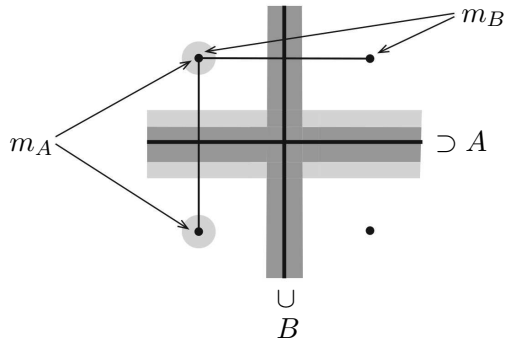}
{\small
\put(-1,99){$m_B$}
\put(-172,50){$m_A$}
\put(-11,50){$\supset A$}
\put(-71,-9){$\cup$}
\put(-71,-21){$B$}
}
\caption{A more detailed illustration of the spheres in figure \ref{spheres} near one of the intersection points. As above, the Clifford torus is shown as $4$ points. $A_n,B_n$ are contained in neighborhoods of the $2$-spheres which are indicated with darker shading. Lighter shading shows the range of a small ambient isotopy.}
\label{isotopy figure}
\end{figure}

This calculation is unchanged by a small isotopy which moves $A_n$ to $\widetilde A_n$. Indeed, the expression for $m_B$, read off in a $3$-dimensional slice, is unchanged (see figure \ref{isotopy figure}). The ambient isotopy applied to the meridian $m_A$ moves it to a curve $\widetilde m_A$ which is contained in a small tubular neighborhood of $m_A$. The expression for $\widetilde m_A$ in the complement of $\widetilde A_n$ is identical to the expression for $m_A$ in the complement of $A_n$, since they are related by an ambient isotopy. Moreover, since $\widetilde m_A$  is contained in a small tubular neighborhood of $m_A$, these two curves are isotopic within this tubular neighborhood. In particular, $m_A$ is isotopic to $\widetilde m_A$ in the complement of $\widetilde A_n$ and of $B_n$. It follows that the expression for $m_A$ as an $l-$fold commutator of the meridians of $A_n$ is unchanged by the isotopy.

Then $[m_A, m_B]$ is a commutator of length $2l$ in ${\pi}_1({\mathbb R}^4\smallsetminus(\widetilde A_n\sqcup B_n))$.    For the remainder of the argument fix an integer $k>2l=N$.
$[m_A, m_B]$ is seen to be non-trivial in ${\pi}_1(Y )/({\pi}_1(Y ))^k\cong  F_N/(F_N)^k$. 
  (The underlying reason for this is that the only relations in $F_N/(F_N)^k$ are: $\{$all commutators of length $k$ are trivial$\}$. A shorter commutator $[m_A, m_B]$, of length $2l<k$ is not a consequence of these relations. A rigorous proof uses the Magnus expansion \cite{M}.) 

This contradiction with (\ref{relation}) completes the proof of Lemma \ref{lemma}.  \qed

Returning to the proof of theorem \ref{theorem}, define the Cantor set $C$ to be the union of $A,B$. Suppose there is an arbitrarily small ambient isotopy, pushing $C$ off of itself. In particular, then $A$ may be isotoped off of $B$ by a small ambient isotopy. This contradiction with Lemma 2 concludes the proof of Theorem 1 in the case $d=4$. 

The construction in higher dimensions $d$ starts from a $(d-3)$-spun Bing decomposition. 
Analogously to the case $d=4$, take two standard $(d-2)$-spheres in general position (intersecting in a $(d-4)$-sphere) in ${\mathbb R}^d$. All steps in the proof above go through: the Clifford torus in figure 2 (corresponding to a codimension$=(d-4)$-slice) gives a relation $[m_A,m_B]=1$ in the fundamental group of the complement. The contradiction with the Stallings theorem relies only on homological information which holds due to Alexander duality. 
\qed

{\em Remark.}  
Michael Freedman suggested a refinement of the construction in theorem \ref{theorem},  giving a {\em locally} sticky Cantor set. That is, the intersection of such a Cantor set with any open set is itself sticky with respect to $\epsilon$-isotopies, for some $\epsilon$ depending on the open set. The idea is to start with the Hopf link and then replace each component with four: a Bing pair, and in addition the two meridional circles, then iterate - always replacing each component with four. 
(In the usual Bing decomposition each component is replaced with a Bing double, shown in figure \ref{Bing figure}. In addition, now the two meridional curves are also included.) This gives a shrinkable decomposition.
The idea is then to use a spun version of it, with ``minimal'' intersections of the $2$-spheres$\, \times D^2$ in $4$-space to define the desired Cantor set.
(The meridional curves create linking numbers, resulting in intersections between the spheres in $4$-space at all scales. )

{\bf Acknowledgments}.  
I would like to thank Bob Edwards, Mike Freedman and Dick Sher for discussions on wild topology. 
I also would like to thank Bob Edwards for bringing the question of existence of sticky Cantor sets to my attention. 

This research was supported in part by NSF grant DMS-1309178 and by the Simons Foundation grant 304272.
I would like to thank the IHES for hospitality and support (NSF grant 1002477).


\begin{thebibliography}{10}

\bibitem{A} L. Antoine, {\em Sur l'homeomorphisme de deux figures et leurs voisinages}, Journal Math Pures et appl., 4
(1921), 221-325.

\bibitem{B} R.H. Bing, {\em A Homeomorphism Between the 3-Sphere and the Sum of Two Solid Horned Spheres},  Ann. Math. 56, (1952), 354-362.

\bibitem{CW} J.W. Cannon and D.G. Wright, 
{\em Slippery Cantor sets in $E^n$}, 
Fund. Math. 106 (1980), 89-98. 

\bibitem{D0} R.J. Daverman, {\em  On the absence of tame disks in certain wild cells},  Geometric topology (Proc. Conf., Park City, Utah, 1974), pp. 142-155. Lecture Notes in Math., Vol. 438, Springer, Berlin, 1975. 

\bibitem{D} R.J. Daverman, Decompositions of manifolds. Pure and Applied Mathematics, 124. Academic Press, Inc., Orlando, FL, 1986. 

\bibitem{D1} R.J. Daverman, {\em Problems about finite-dimensional manifolds}, Open problems in topology, 431-455, North-Holland, Amsterdam, 1990.

\bibitem{DE} R.J. Daverman and R.D.  Edwards, 
{\em Wild Cantor sets as approximations to codimension two manifolds},
Topology Appl. 26 (1987),  207-218. 

\bibitem{DV} R.J. Daverman and G.A. Venema, Embeddings in manifolds. 
Graduate Studies in Mathematics, 106. American Mathematical Society, Providence, RI, 2009.

\bibitem{Edwards} R.D. Edwards, {\em Suspensions of homology spheres}, ArXiv 0610573.

\bibitem{L} L.L. Lininger, {\em Actions on $S^n$}, Topology 9 (1970), 301-308.

\bibitem{M} J. Milnor, {\em Link Groups}, Ann. Math. 59 (1954), 177-195.




\bibitem{S} J. Stallings, {\em Homology and central series of groups}, J. Algebra 2 (1965) 170-181. 

\bibitem{W} D.G. Wright,  {\em Pushing a Cantor set off itself},
Houston J. Math. 2 (1976),  439-447.


\end{thebibliography}
\end{document}